\documentclass[letterpaper,runningheads]{llncs}

\usepackage{amssymb}
\usepackage{graphicx}
\usepackage{color}

\usepackage{url}
\urldef{\mailsa}\path|{d.andreagiovanni,gleixner}@zib.de|
\newcommand{\keywords}[1]{\par\addvspace\baselineskip
\noindent\keywordname\enspace\ignorespaces#1}

\usepackage{amsmath}

\usepackage{tabularx,booktabs,lscape,array}

\usepackage{xspace}

\usepackage{url}
\usepackage[bookmarks,implicit=true,colorlinks,linkcolor=black,citecolor=black,urlcolor=blue]{hyperref}
\hypersetup{%
pdftitle=Towards an accurate solution of wireless network design problems,
pdfauthor={Fabio D'Andreagiovanni and Ambros M.~Gleixner},
pdfkeywords={Linear Programming, Precise Solutions, Network Design, Wireless Telecommunications Systems}
}

\newcommand{\feasible}{$\checkmark$}
\newcommand{\infeasible}{$\varnothing$}

\newcommand{\nbsc}[1]{\scalebox{0.9}{#1}\xspace}

\newcommand{\e}[1]{\ensuremath{\!\cdot\!10^{#1}}}

\newcommand{\T}{\mathsf{T}\hspace*{-0.2ex}}

\renewcommand{\leq}{\leqslant}
\renewcommand{\geq}{\geqslant}

\newcommand{\MIP}{\nbsc{MIP}}
\newcommand{\SIR}{\nbsc{SIR}}
\newcommand{\MINLP}{\nbsc{MINLP}}

\newcommand{\LP}{\nbsc{LP}}
\newcommand{\LU}{\nbsc{LU}}

\newcommand{\NETLIB}{\nbsc{NETLIB}}

\newcommand{\CPLEX}{\nbsc{CPLEX}}
\newcommand{\SCIP}{\nbsc{SCIP}}
\newcommand{\SOPLEX}{\nbsc{SoPlex}}
\newcommand{\QSOPTEX}{\nbsc{QSopt\_ex}}
\newcommand{\GMP}{\nbsc{GMP}}
\newcommand{\EGLIB}{\nbsc{EGlib}}

\definecolor{darkred}{rgb}{0.8,0,0}

\newcolumntype {Y}            {>{\raggedright\arraybackslash}X}
\newcolumntype {Z}            {>{\raggedleft\arraybackslash}X}
\newcolumntype {F} [1]        {D{.}{.}{#1}}
\newcolumntype {N} [1]        {D{.}{.}{#1.0}<{\hspace{-0.8ex}}}
\newcolumntype{M}{>{\raggedleft\arraybackslash}X}
\begin{document}

\mainmatter  

\title{Towards an accurate solution of wireless network design problems%
  \thanks{This is the authors' final version of the paper published in: R. Cerulli, S. Fujishige, A.R. Mahjoub (Eds.),  Combinatorial Optimization - 4th International Symposium, ISCO 2016, LNCS 9849, pp. 135-147, 2016, DOI: 10.1007/978-3-319-45587-7\_12.
The final publication is available at Springer via http://dx.doi.org/10.1007/978-3-319-45587-7\_12.
This work was partly conducted within the Research
Campus Modal funded by the German Federal Ministry of Education and Research (Grant no. 05M14ZAM). It was also partially supported
by the \emph{Einstein Center for Mathematics Berlin} (ECMath) through Project MI4 (ROUAN) and by the \emph{German Federal Ministry of Education and Research} (BMBF) through Project VINO (Grant no. 05M13ZAC) and Project \emph{ROBUKOM} (Grant no. 05M10ZAA) \cite{BaEtAl14}.}%
}

\titlerunning{Accurate solution of WND problems}

%
%

\author{Fabio D'Andreagiovanni\inst{1,}\inst{2,}\inst{3} \and Ambros M.\ Gleixner\inst{1}%
}
\authorrunning{F. D'Andreagiovanni \and A.M.\ Gleixner}

\institute{
Dept. of Mathematical Optimization, Zuse Institute Berlin (ZIB)\\
Takustr. 7, 14195 Berlin, Germany\\
\and
DFG Research Center \textsc{Matheon}, Einstein Center for Mathematics (ECMath)\\
Stra{\ss}e des 17.\ Juni 136, 10623 Berlin, Germany\\
\and
Institute for System Analysis and Computer Science, National Research Council of Italy (IASI-CNR), via dei Taurini 19, 00185 Roma, Italy\\
\mailsa
}

%
%

\toctitle{Lecture Notes in Computer Science}
\tocauthor{Authors' Instructions}
\maketitle

\begin{abstract}
The optimal design of wireless networks has been widely studied in the literature and many optimization models have been proposed over the years.
However, most models directly include the signal-to-interference ratios representing service coverage conditions.
This leads to mixed-integer linear programs with constraint matrices containing tiny coefficients that vary widely in their order of magnitude.
These formulations are known to be challenging even for
state-of-the-art solvers: the standard numerical precision supported by these solvers is
usually not sufficient to reliably guarantee feasible solutions. Service coverage errors are thus commonly present.
Though these numerical issues are known and become evident even for small-sized instances, just a very limited number of papers has tried to tackle them, by mainly investigating alternative non-compact formulations in which the sources of numerical instabilities are eliminated.
In this work, we explore a new approach by investigating how recent advances in exact solution algorithms for linear and mixed-integer programs over the rational numbers
can be applied to analyze and tackle the numerical difficulties arising in wireless network design models.

\end{abstract}


\section{Introduction}
\label{sec:intro}

In the last decade, the presence of wireless communications in our everyday life has greatly expanded and wireless
networks have thus increased in number, size and technological complexity.
In this context, the traditional design approach adopted by professionals, based on trial-and-error supported by
simulation, has exhibited many limitations. This approach is in particular not able to pursue an efficient exploitation
of scarce and precious radio resources, such as frequency channels and channel bandwidth, and the need for exact
mathematical optimization approaches has increased.

The problem of designing a wireless network can be essentially described as that of configuring a set of transmitters in
order to cover with a telecommunication service a set of receivers, while guaranteeing a minimum quality of service.
Over the years, many optimization models have been proposed for designing wireless networks (see
\cite{DA12,DAMaSa13,KeOlRa10} for an introduction).
However, most models have opted for so-called \emph{natural formulations}, which directly include the formulas used to
assess service coverage conditions. This leads to the definition of mixed-integer programs whose constraint matrices
contain tiny coefficients that  greatly vary in their order of magnitude. Furthermore, the natural formulations commonly
include also the notorious big-M coefficients to represent disjunctive service coverage constraints.  These formulations are known to be
challenging even for state-of-the-art solvers. Additionally, the standard numerical precision supported by these solvers
is usually not sufficient to reliably guarantee feasible solutions \cite{MaRoSm07}. If returned solutions are verified
in a post-optimization phase, it is thus common to find service coverage errors.

Though these numerical issues are known and can be found even in the case of instances of small size, it is interesting to note that just a very limited number of papers has tried to tackle them. In particular, the majority of these works rely on the definition of alternative non-compact formulations that are able to reduce the numerical drawbacks of natural formulations (see the next section for a review of the main approaches).

In contrast to these works, we propose here a new approach: we investigate how
recent advances in exact solution algorithms for (integer) linear programs over the rational numbers
can be applied to analyze and tackle the numerical difficulties arising in wireless network design.


\smallskip
\noindent
Our main original contributions are in particular:
\begin{enumerate}
    \item we present the first formal discussion about why even effective state-of-the-art solvers fail to correctly
          discriminate between feasible and infeasible solutions in wireless network design;
    \item we assess, for the first time in literature, both formally and computationally the actual benefits coming from
          scaling the very small coefficients involved in natural formulations; coefficient scaling is a practice that is adopted
          by many professionals and scholars dealing with wireless network design, with the belief of eliminating numerical
          errors;
        we show that just adopting scaling is not
sufficient to guarantee accurate feasibility of solutions returned by floating-point solvers;
    \item we show how extended-precision solvers can be adopted to check the correctness of solutions returned by floating-point
solvers and, if errors are present, to get correct valorization of the continuous variables of the
problem.
\end{enumerate}
Our computational experiments are made over a set of realistic instances defined in collaboration with a major
European telecommunication company.


\medskip
\noindent
The remainder of this paper is organized as follows:
in Section~\ref{sec:WND}, we formally characterize the wireless network design problem and introduce the natural
formulations; in Section~\ref{sec:algorithms}, we discuss the question of accuracy in Mixed Integer Programming (\MIP) solvers,
addressing in particular the issues arising in wireless network design; in Section~\ref{sec:experiments}, we present our
computational experiments over realistic network instances.


\section{The Wireless Network Design Problem}
\label{sec:WND}

For modeling purposes, a wireless network can be described as a set of
transmitters $T$ that provide a telecommunication service to a set of
receivers $R$.  Transmitters and receivers are
characterized by a location and a number of radio-electrical
parameters (e.g., power emission and transmission frequency).  The {\em Wireless
Network Design Problem}  (WND)
consists in establishing the location and suitable values for the parameters of
the transmitters with
the goal of optimizing an objective function that expresses the interest of the
decision maker: common objectives are the maximization of a revenue function
associated with wireless service coverage or, assuming a green-network
perspective, the minimization of the total power emission of the network
transmitters. For an exhaustive introduction to the WND, we refer the reader to
\cite{DA12,DAMaSa13,KeOlRa10}.

Given a receiver $r \in R$ that we want to cover with service, we must choose a single transmitter $s \in S$, called
\emph{server}, that provides the telecommunication service to $r$. Once the server of a receiver is chosen, all the
other transmitters are \emph{interferers} and deteriorate the quality of service obtained by $r$ from its server~$s$.

From an analytical point of view, if we denote by $p_t$ the power emission of a transmitter $t \in T$, a receiver $r \in R$ is considered covered with service (or briefly
\emph{served}) when the ratio of the service power to the sum of the interfering
powers
(\emph{Signal-to-Interference Ratio} - \emph{SIR}) is above a threshold $\delta >
0$, that depends on the desired quality of service \cite{Ra01}:
%
\begin{equation}
\label{eq:firstSIRineq}
SIR_{r s}(p) = \frac{a_{r s(r)} \cdot p_{s(r)}}
{N + \sum_{t \in T\setminus\{s(r)\}} a_{rt} \cdot p_t}
\hspace{0.1cm}
\geq
\hspace{0.1cm}
\delta \;  .
\end{equation}

\noindent
In this inequality:
i) $s(r) \in T$ is the server of receiver $r$;
ii) the power $P_{t}(r)$ that $r$ receives from a transmitter $t \in T$ is
proportional to the emitted power $p_t$ by a factor $a_{rt} \in [0,1]$, i.e.
$P_{t}(r) = a_{rt}\cdot p_t$. The factor $a_{rt}$ is called {\em fading
coefficient} and summarizes the reduction in power that a signal experiences
while propagating from $t$ to $r$ \cite{Ra01};
iii) in the denominator, we highlight the presence of the system noise  $N >
0$ among the interfering signals.
%

\noindent
By simple algebra operations, inequality (\ref{eq:firstSIRineq}) can be
transformed into the
following linear inequality, commonly called \emph{SIR inequality}:
\begin{equation}\label{eq:secondSIRineq}
a_{rs(r)}  \cdot p_{s(r)} - \delta \sum_{t \in T \setminus\{s(r)\}} a_{rt}  \cdot p_t
\hspace{0.1cm}
\geq
\hspace{0.1cm}
\delta \cdot N \;  .
\end{equation}

\noindent
Since service coverage assessment is a central element in the design of any
wireless network,  the SIR inequality constitutes the core of any optimization
problem used in wireless network design. If we just focus attention on setting
power emissions, we can define the so-called \emph{Power Assignment Problem}
(PAP), in which we want to fix the power emission of each transmitter in order to serve a set of receivers, while minimizing the sum of all power emissions. By introducing a non-negative
decision variable $p_t \in [0,P_{\max}]$ to represent the feasible power
emission range of a transmitter $t \in T$, the PAP can be easily formulated as
the following pure Linear Program (\LP):
\begin{align}
\min
&
\sum_{t \in T} p_t
&&
\tag{PAP}
\nonumber
\\
&
a_{r s(r)}  \cdot p_{s(r)} - \delta \sum_{t \in T \setminus\{s(r)\}} a_{rt}  \cdot
p_t
\hspace{0.1cm}
\geq
\hspace{0.1cm}
\delta \cdot N
&&
\forall \hspace{0.1cm} r \in R
\label{PAPsir}
\\
&
0 \leq p_t \leq P_{\max}
&&
\forall \hspace{0.1cm}  t \in T \;  ,
\end{align}

\noindent
where \eqref{PAPsir} are the SIR inequalities associated with receivers to
be served.

In a hierarchy of WND problems (see \cite{DA12,MaRoSm07} for details), the PAP constitutes a basic WND problem that
lies at the core of virtually all more general WND problems.
A particularly important generalization of the PAP is constituted by the \emph{ Scheduling and
Power Assignment Problem (SPAP)} \cite{DA12,DAMaSa13,MaRoSm07,MaRoSm06}, where, besides the power emissions,
it is also necessary to choose the assignment of a served receiver to a transmitter in the network that acts as server of the receiver.
This can be easily modeled by introducing 0-1 service assignment variables, obtaining the following natural formulation:
\begin{align}
\max\;\;
&
\sum_{r \in R} \sum_{t \in T} \pi_{t} \cdot x_{rt}
&&
\tag{SPAP}
\\
&
a_{r s}  \cdot p_{s} - \delta  \sum_{t \in T \setminus\{s\}} a_{rt}  \cdot  p_t
+ M  \cdot  (1 - x_{rs})
\hspace{0.1cm}
\geq
\hspace{0.1cm}
\delta \cdot N
&&
\forall \hspace{0.1cm} r \in R, s \in T
\label{SPAPsir}
\\
&
\sum_{t \in T} x_{rt} \leq 1
&&
\forall \hspace{0.1cm} r \in R
\label{SPAPgub}
\\
&
0 \leq p_t \leq P_{\max}
&&
\forall \hspace{0.1cm} t \in T
\\
&
x_{rt} \in \{0,1\}
&&
\forall \hspace{0.1cm} r \in R, t \in T \;  ,
\end{align}

\noindent
which includes:
i) additional binary variables $x_{rt}$ to represent that receiver $r$ is served by transmitter $t$;
ii) modified SIR inequalities, defined for each possible server transmitter $s \in T$ of a receiver r, including large constant values $M > 0$ to activate/deactivate the corresponding SIR inequalities (as expressed by the constraint \eqref{SPAPgub} each user may be served by at most one transmitter and thus at most one SIR inequality must be satisfied for each receiver);
iii) a modified objective function aiming at maximizing the revenue obtained from serving transmitters (every receiver grants a revenue $\pi_{t} > 0$).
%
%


\vspace*{-1ex}

\paragraph*{\bf Drawbacks of SIR-based formulations.}

The natural (mixed-integer) linear programming formulations associated with the PAP and the SPAP and based on the direct
inclusion of the SIR inequalities are widely adopted for the WND in different application contexts, such as DVB-T, (e.g.,
\cite{MaRoSm07,MaRoSm06}), UMTS (e.g.,
\cite{AmEtAl06}),
WiMAX (e.g., \cite{DA12,DAMaSa13}).
In principle, such formulations can be solved by MIP solvers, but, as clearly pointed out in works like
\cite{DA12,DAMaSa13,KeOlRa10,MaRoSm07}, in practice:
\begin{itemize}
  \item the fading coefficients may vary in a wide range leading to very ill-conditioned coefficient matrices (for example, in the case of DVB instances, difference between coefficients may exceed 90 decibels) that make the solution process \emph{numerically unstable};
  \item in the case of SPAP-like formulations, the big-{\em M} coefficients lead to extremely weak bounds that may greatly decrease the effectiveness of solvers implementing state-of-the-art versions of branch-and-bound techniques;
  \item the resulting coverage plans are often unreliable and may contain errors, i.e. SIR constraints recognized as satisfied by an MIP solver actually reveal to be violated.
\end{itemize}

\noindent
Though these issues are known, it is interesting to note that just a limited number of works in the wide literature about WND has tried to tackle them and natural formulations are still widely used. We refer the reader to \cite{DA12,KeOlRa10} for a review of works that have tried to tackle these drawbacks and we recall here some more relevant ones.
One of the first works that has identified the presence and effects of numerical issues in WND is \cite{MaRoSm06}, where a GRASP algorithm is proposed to solve very large instances of the SPAP, arising in the design of  DVB-T networks. Other exact solution approaches have aimed at eliminating the source of numerical instabilities (i.e., the fading and big-M coefficients) by considering non-compact formulations: in \cite{CaEtAl11}, a formulation based on cover inequalities is introduced for a maximum link activation problem; in \cite{DA12,DAMaSa13}, it is instead shown how using a \emph{power-indexed} formulation, modeling power emissions by discrete power variables allows to define a peculiar family of generalized upper bound cover inequalities that provide (strong) formulations. In \cite{DA12}, it is also presented an alternative formulation based on binary expansion of variables, which can become strong in some relevant practical cases, thanks to the superincreasing property of the adopted expansion coefficients.
In \cite{DAMaSa11}, it is proposed the definition of a non-compact formulation purely based on assignment variables that relates to a maximum feasible subsystem problem. Finally, in \cite{DA11}, the numerical instabilities are addressed by the definition of a genetic heuristic exploiting the discretization of power emissions.

According to a widespread belief,
numerical instabilities in WND may be eliminated by multiplying all the fading coefficients of the problem by a large power of 10 (typically $10^{12}$).
However, in our direct experience with real-world instances of several wireless
technologies (e.g., DVB-T \cite{DAMaSa13}, WiMAX \cite{DA12,DAMaSa13}), this did neither improve the performance of the solver nor of the quality of solutions found, which were still subject to coverage errors.


\section{Numerical accuracy in linear and mixed-integer linear programming solvers}
\label{sec:algorithms}




Wireless network design problems are not only combinatorially complex, but as was argued before, also numerically
sensitive.
%
State-of-the-art \MIP solvers employ floating-point arithmetic, hence their arithmetic computations are subject to
round-off errors.  This makes it necessary to allow for small violations of the constraints, bounds,
and integrality requirements when checking solutions for feasibility.
To this end, \MIP solvers typically use a combination of absolute and relative tolerance to
define their understanding of feasibility.
A linear inequality $\alpha^\T x \leq \alpha_0$ is considered as satisfied by a point $x^*$ if
\begin{align}
  \label{equ:feasibility}
  \frac{\alpha^\T x^* - \alpha_0}{ \max\{ |\alpha^\T x^*|, |\alpha_0|, 1 \} } \leq
\epsilon_{\text{feas}}
\end{align}
with a feasibility tolerance $\epsilon_{\text{feas}} > 0$.\footnote{%
This is the definition of feasibility used by the academic \MINLP solver
\SCIP~\cite{Achterberg2009,scip}.
While we do not know for certain the numerical definitions used by
closed-source commercial solvers, we think
that they follow a similar practice.}
If the activity $\alpha^\T x^*$ and right-hand side $\alpha_0$ are below one in absolute
value, an absolute violation of up
to~$\epsilon_{\text{feas}} > 0$ is allowed.  Otherwise, a relative tolerance is
applied and larger violations are
accepted.  Typically, $\epsilon_{\text{feas}}$ ranges between $10^{-6}$ and
$10^{-9}$.

\vspace*{-1ex}

\paragraph*{\bf Feasibility of SIR inequalities.}

When employing floating-point arithmetic to optimize wireless network design problems containing SIR inequalities,
care is required when enforcing and checking their feasibility.
First, since the coefficients and right-hand side of the linearized SIR inequality~\eqref{eq:secondSIRineq} are
significantly below $10^{-9}$ in absolute value, the inequality \eqref{equ:feasibility} results in a very loose definition of
feasibility.  The allowed absolute violation may be larger than the actual right-hand side.

Second, though the original SIR inequality~\eqref{eq:firstSIRineq} is equivalent to its linear
reformulation~\eqref{eq:secondSIRineq}, if we check their violation with respect to numerical tolerances,
they behave differently.  Indeed, an (absolute) violation~$\epsilon_{\text{linear}}
= \delta N - (a_{rs} p_{s} - \delta \sum_{t \in T \setminus\{s\}} a_{rt} p_t)$ of~\eqref{eq:secondSIRineq}
corresponds to a much larger violation of~\eqref{eq:firstSIRineq}, since
\begin{align}
  \epsilon_{\text{SIR}}
  = \delta - \frac{a_{r s} p_{s}}{N + \sum_{t \in T\setminus\{s\}} a_{rt} p_t}
  = \frac{\epsilon_{\text{linear}}}{N + \sum_{t \in T\setminus\{s\}} a_{rt} p_t}
\end{align}
and the sum of noise and interference signal $N + \sum_{t \in T\setminus\{s\}} a_{rt} p_t$ typically has an order of
$10^{-9}$ or smaller.  In combination with the feasibility tolerances promised by standard \MIP solvers ($\approx10^{-9}$),
this would at best guarantee violations in the order of $1$ for the original problem formulation.

\vspace*{-1ex}

\paragraph*{\bf The impact of scaling.}
\label{sec:scaling}

Internally, MIP solvers may apply scaling factors to rows and columns of the
constraint matrix in order to improve the numerical stability.  Primarily, this aims at improving the
condition numbers of basis matrices during the solution of \LP{}s.

However, from~\eqref{equ:feasibility} it becomes apparent that an external, a
priori scaling of constraints by the user can change the very definition
of feasibility: if the activity and right-hand side are significantly below 1
in absolute value, then scaling up tightens the feasible region.  Precisely,
with a scaling factor $S > 1$, if $|S\alpha^\T x^*| < 1$ and $|S\alpha_0| < 1$, then
\begin{align}
  \frac{S(\alpha^\T x^* - \alpha_0)}{ \max\{ |S\alpha^\T x^*|, |S\alpha_0|, 1 \} } \leq
\epsilon_{\text{feas}}
  \Leftrightarrow
  \frac{(\alpha^\T x^* - \alpha_0)}{ \max\{ |\alpha^\T x^*|, |\alpha_0|, 1 \} } \leq
\frac{\epsilon_{\text{feas}}}{S},
\end{align}
and the absolute tolerance can be decreased by a factor of~$1/S$.
This can then be
used to arrive at a sufficiently strict definition of feasibility for
constraints with very small coefficients, such as the SIR
inequalities~\eqref{eq:secondSIRineq}.

\vspace*{-1ex}

\paragraph*{\bf Advances in exact \LP and \MIP solving.}
\label{sec:exactmethods}

Although the floating-point numerics used in today's state-of-the-art MIP
solvers yield reliable results for the majority
of problems and applications, there are cases in which results of higher
accuracy are desired or needed, such as
verification problems, computer proofs, or simply numerically instable
instances.
In the following we will review recent advances in methods for solving \LP{s}
and \MIP{}s exactly over the rational numbers.

Trivially, of course, one can obtain an exact solution algorithm by performing
all computations in exact arithmetic.
However, for all but a few instances of interest, this idea is not sufficiently
performant.
As a starting point, it has been observed that \LP bases returned by
floating-point solvers are often
optimal for real world problems \cite{Dhiflaoui2003}.  For example, \cite{Koch2004} could
compute optimal bases to all of the \NETLIB \LP instances 
using only floating-point LP solvers and subsequently certifying them in exact
rational arithmetic.

Following these observations,
Applegate et al. \cite{Applegate2007b} developed a
simplex-based general-purpose exact \LP solver,
\QSOPTEX, that exploits this behavior to achieve fast computation times on
average. If an optimal basis is not
identified by the double-precision subroutines, more simplex pivots are
performed using increased levels of precision until the exact rational
solution is identified.  For more details, see~\cite{Espinoza2006}.

Recently, Gleixner et al.~\cite{GleSteWol2012,GleSteWol2016} have developed an
iterative refinement procedure for solving \LP{}s
with high accuracy by solving a sequence of closely related \LP{s} in order to compute
primal and dual correction terms.
The procedure avoids rational \LU factorizations and \LP solves in extended precision and
hence often computes solutions with only tiny violations
faster than \QSOPTEX.
Although not an exact method in itself, it can be used to speed up \QSOPTEX
significantly.

Finally, exact \LP solving is a crucial subroutine for solving \MIP{}s exactly.  Once a promising assignment for the
integer variables has been found, an exact \LP solver can be used to compute feasible values for the continuous
variables or prove that this integer assignment does not admit a fully feasible solution vector.

The majority of \LP{}s within a \MIP solution process, however, is solved to bound the objective value of the optimal
solution.  Solving these exactly does provide safe dual bounds, but can result in a large slow-down.
The key to obtain a faster exact \MIP solver
is to avoid exact \LP solving by correcting the dual solution obtained from a floating-point \LP solver,
see~\cite{neumaier2004safe}.  Cook et al.~\cite{Cook2011,CookKochSteffyWolter2013} have followed
this approach to develop an exact branch-and-bound algorithm available as an extension of the solver \SCIP~\cite{scip}.

In the following section, we will investigate empirically how these tools can be applied to analyze and
address the numerical difficulties encountered in solving wireless network design problems.


\section{Computational experiments}
\label{sec:experiments}

The goal of our experiments was twofold: first, in order to test whether \MIP solvers can be reliably
used as decision tools for wireless network design models as introduced in Sec.~\ref{sec:WND},
we analyzed the accuracy of primal solutions returned by a state-of-the-art \MIP solver;
second, we investigated the practical applicability and performance of the exact solution
methods described in Sec.~\ref{sec:exactmethods}.

\vspace*{-1ex}

\paragraph*{\bf Experimental setup.}
The experiments were conducted on a computer with a 64bit Intel Xeon E3-1290 v2 CPU (4 cores, 8 threads)
at 3.7\,GHz with 8\,MB cache and 16\,GB main memory.
We ran all jobs separately to avoid random noise in the measured running time that might be caused by cache-misses if
multiple processes share common resources.  We used \CPLEX~12.5.0.0~\cite{cplex} (default, deterministic parallel with up to four threads),
\QSOPTEX~2.5.10~\cite{Applegate2007b} with \EGLIB~2.6.10 and \GMP~4.3.1~\cite{gmp}, and \SOPLEX~2.0~\cite{soplex}
with \GMP~5.0.5~\cite{gmp} (both single-thread).

\vspace*{-1ex}

\paragraph*{\bf Test instances.}
We performed our experiments on realistic instances of a WiMAX network,
defined in cooperation with a major European telecommunications company.
The instances correspond to various scenarios of a single-frequency network
adopting a single transmission scheme.%
\footnote{For more details on WiMAX networks, see~\cite{DA12}.}
For each instance, we solved the corresponding SPAP model from Sec.~\ref{sec:WND}.

We considered ten instances with between 100 and 900~receivers ($|R|$) and
between 8 and 45~transmitters ($|T|$).  The maximum emission
power of each transmitter ($P_{\max}$) was set equal to 30\,dBmW and the SIR threshold ($\delta$)
was between 8\,dB and 11\,dB.%
\footnote{The smallest \MIP has 808~variables, 900~constraints, and 8\,000~nonzeros, the largest instance contains 32\,436~variables, 33\,300~constraints, and 1\,231\,200~nonzeros.}

\vspace*{-1ex}

\paragraph*{\bf Accuracy of \MIP solutions.}
In our first experiment, we ran \CPLEX with a time limit of one hour (because of the combinatorial complexity of
the problems, only the smallest instances can be solved to optimality within this limit) and checked the feasibility
of the best primal solution returned.
Table~\ref{tab:fixing0} shows the results for the unscaled instances, Table~\ref{tab:fixing12} shows the
results for the instances with the linearized \SIR inequalities~\eqref{eq:secondSIRineq} multiplied by $S=10^{12}$
as in Sec.~\ref{sec:scaling}.

\begin{table}[t]
  \centering
  \caption[ ]{A posteriori check and exact verification of binary assignments from floating-point MIP solutions for
    instances without scaling.}
   \setlength{\tabcolsep}{5pt}
  \scriptsize
  \begin{tabularx}{\textwidth}{rrrrrlrrrcr}
    \toprule
    \multicolumn{4}{c}{instance} & & \multicolumn{4}{c}{post processing} & \multicolumn{2}{c}{exact \LP} \\
    \cmidrule(lr){1-4}\cmidrule(lr){6-9}\cmidrule(lr){10-11}
    $|R|$ & $|T|$ & $\alpha_{\min}$ & $\alpha_{\max}$  & obj. & linear viol. & \SIR viol. & served & \!\!\!unserved & stat. & time \\
    \midrule
      100 &    8 & $4\e{-17}$ & $4\e{-8}$ &  41 & $1.7\e{-10}$ &$12.6$ &  13 &  28 & \infeasible &   0.2 \\
  169 &   20 & $1\e{-19}$ & $3\e{-8}$ &  73 & $6.4\e{-11}$ & $6.3$ &   1 &  72 & \infeasible &  20.3 \\
  225 &   20 & $2\e{-19}$ & $2\e{-8}$ & 176 & $1.2\e{-10}$ & $6.3$ &   5 & 171 & \infeasible &  10.0 \\
  256 &   40 & $4\e{-19}$ & $3\e{-8}$ & 155 & $9.0\e{-11}$ & $6.3$ &  15 & 140 & \infeasible & 103.0 \\
  400 &   25 & $8\e{-20}$ & $2\e{-8}$ & 373 & $1.2\e{-10}$ & $6.3$ &   7 & 366 & \infeasible &  55.7 \\
  400 &   40 & $8\e{-20}$ & $2\e{-8}$ & 301 & $9.0\e{-11}$ & $6.3$ &  13 & 288 & \infeasible & 233.7 \\
  441 &   45 & $8\e{-20}$ & $2\e{-8}$ & 312 & $1.0\e{-10}$ & $6.3$ &  15 & 297 & \infeasible & 440.5 \\
  529 &   40 & $8\e{-20}$ & $2\e{-8}$ & 421 & $9.0\e{-11}$ & $6.3$ &  13 & 408 & \infeasible & 337.1 \\
  625 &   25 & $2\e{-17}$ & $5\e{-5}$ & 280 &  $1.9\e{-9}$ & $6.0$ & 225 &  55 & \infeasible & 113.2 \\
  900 &   36 & $2\e{-20}$ & $9\e{-9}$ & 890 & $7.7\e{-11}$ & $2.5$ &  14 & 876 & \infeasible & 660.1 \\
    \bottomrule
  \end{tabularx}
  \label{tab:fixing0}
\end{table}

\begin{table}[t]
  \centering
  \caption[ ]{A posteriori check and exact verification of binary assignments from floating-point MIP solutions for
    instances scaled with $10^{12}$.}
   \setlength{\tabcolsep}{5pt}
  \scriptsize
  \begin{tabularx}{\textwidth}{rrrrrlrrrcr}
    \toprule
    \multicolumn{4}{c}{instance} & & \multicolumn{4}{c}{post processing} & \multicolumn{2}{c}{exact \LP} \\
    \cmidrule(lr){1-4}\cmidrule(lr){6-9}\cmidrule(lr){10-11}
    $|R|$ & $|T|$ & $\alpha_{\min}$ & $\alpha_{\max}$  & obj.\;\; & linear viol. & \SIR viol. & served & unserved & \;stat. & time \\
    \midrule
      100 &    8 & $4\e{-5}$ & $4\e{5}$ &  28 & $7.1\e{-17}$ & $4.6\e{-6}$ &  24 &  4 & \feasible &   0.1 \\
  169 &   20 & $1\e{-7}$ & $3\e{5}$ &  44 & $7.3\e{-17}$ & $8.0\e{-6}$ &  43 &  1 & \feasible &   2.1 \\
  225 &   20 & $3\e{-7}$ & $2\e{5}$ &  42 & $6.2\e{-17}$ & $7.0\e{-6}$ &  38 &  4 & \feasible &   0.9 \\
  256 &   40 & $4\e{-7}$ & $3\e{5}$ &  72 & $8.1\e{-17}$ & $5.1\e{-6}$ &  62 & 10 & \feasible &  11.7 \\
  400 &   25 & $8\e{-8}$ & $2\e{5}$ &  77 & $6.7\e{-17}$ & $1.4\e{-5}$ &  71 &  6 & \feasible &   5.5 \\
  400 &   40 & $8\e{-8}$ & $2\e{5}$ &  95 & $6.7\e{-17}$ & $7.9\e{-6}$ &  85 & 10 & \feasible &  20.2 \\
  441 &   45 & $8\e{-8}$ & $2\e{5}$ & 101 & $8.9\e{-16}$ & $4.8\e{-5}$ &  89 & 12 & \feasible &  35.5 \\
  529 &   40 & $8\e{-8}$ & $2\e{5}$ & 101 & $8.8\e{-15}$ & $6.1\e{-4}$ &  96 &  5 & \feasible &  29.9 \\
  625 &   25 & $8\e{-5}$ & $5\e{7}$ & 417 & $1.9\e{-14}$ & $1.9\e{-3}$ & 415 &  2 & \feasible &   6.0\\
  900 &   36 & $2\e{-8}$ & $9\e{4}$ & 202 & $8.1\e{-18}$ & $1.4\e{-6}$ & 200 &  2 & \feasible &  58.0 \\
    \bottomrule
  \end{tabularx}
  \label{tab:fixing12}
\end{table}

The first two columns give the size of each instance, while the second two columns state the smallest and largest absolute
value in the coefficients and right-hand sides of the SIR constraints.  These values differ by up
to~$10^{12}$, a first indicator of numerical instability.  Column ``obj.'' gives the objective value of the
solution at the end of the solving process that we checked, i.e., the number of receivers served by one transmitter.
We report both the maximum violation of the original \SIR
inequalities~\eqref{eq:firstSIRineq} in column ``\SIR viol.'' and their linearization~\eqref{eq:secondSIRineq}.
Both for scaled and unscaled models, the results show that they differ by a factor of up to $10^{12}$.  This demonstrates
that the linearized \SIR inequalities must be satisfied with a very tight tolerance if we want to guarantee a reasonably
small tolerance, $10^{-6}$, say, for the original problem statement.

As it can be seen, the results for the unscaled models
are significantly worse in this respect: although the violation of the linearized constraint looks quite small, the
original \SIR inequalities are strongly violated.  As a result, these solutions cannot be implemented in practice.%
\footnote{Although with this kind of unreliability, this does not matter anymore, note that the numerical difficulties during the
solving process are also reflected in the lower objective values obtained by the unscaled models.}

The column ``served'' states the number of receivers~$r$ served by a transmitter~$s$ for which the corresponding quantity $SIR_{rs}(p)$
is at least $\delta - 10^{-6}$.  This gives the (cardinality of the) subset of receivers that can reliably be served by the
power vector $p$ of the \MIP solution.  It is evident these values are significantly below the claimed objective value of
the \MIP solution for the unscaled models.  Although the situation is much better for the scaled models, also these exhibit
a notable number of receivers that are incorrectly claimed to be served.

\begin{table}[t]
  \centering
  \caption[ ]{Exact computation of the power vector via \QSOPTEX versus iterative refinement via \SOPLEX to a tolerance of $10^{-25}$
    for instances scaled with $10^{12}$.}
   \setlength{\tabcolsep}{6pt}
  \scriptsize
  \begin{tabularx}{38.5em}{rrrrrcrrrr}
    \toprule
    \multicolumn{4}{c}{instance} & & \multicolumn{2}{c}{\QSOPTEX} & \multicolumn{3}{c}{\SOPLEX} \\
    \cmidrule(lr){1-4}\cmidrule(lr){6-7}\cmidrule(lr){8-10}
    $|R|$ & $|T|$ &  $\alpha_{\min}$ & $\alpha_{\max}$ & obj. & stat. & time  & max.\ viol.\ & time & rel. [\%] \\
    \midrule
      100 &    8 & $4\e{-5}$ & $4\e{5}$ &  28 & \feasible &    0.1 & $3.7\e{-29}$ &   0.1 &  $-0.0$ \\
  169 &   20 & $1\e{-7}$ & $3\e{5}$ &  44 & \feasible &    2.1 & $2.4\e{-39}$ &   1.0 & $-52.4$ \\
  225 &   20 & $3\e{-7}$ & $2\e{5}$ &  42 & \feasible &    0.9 & $2.2\e{-29}$ &   0.5 & $-44.4$ \\
  256 &   40 & $4\e{-7}$ & $3\e{5}$ &  72 & \feasible &   11.7 & $1.6\e{-36}$ &   2.5 & $-78.6$ \\
  400 &   25 & $8\e{-8}$ & $2\e{5}$ &  77 & \feasible &    5.5 & $1.8\e{-27}$ &   1.3 & $-76.4$ \\
  400 &   40 & $8\e{-8}$ & $2\e{5}$ &  95 & \feasible &   20.2 & $6.9\e{-40}$ &   4.4 & $-78.2$ \\
  441 &   45 & $8\e{-8}$ & $2\e{5}$ & 101 & \feasible &   35.5 & $3.1\e{-40}$ &   6.3 & $-82.2$ \\
  529 &   40 & $8\e{-8}$ & $2\e{5}$ & 101 & \feasible &   29.9 & $5.7\e{-27}$ &   4.9 & $-83.6$ \\
  625 &   25 & $8\e{-5}$ & $5\e{7}$ & 417 & \feasible &    6.0 & $3.0\e{-29}$ &   2.8 & $-53.3$ \\
  900 &   36 & $2\e{-8}$ & $9\e{4}$ & 202 & \feasible &   58.0 & $5.4\e{-40}$ &  11.7 & $-79.8$ \\

%
%
    \bottomrule
  \end{tabularx}
  \label{tab:ir}
\end{table}

\vspace*{-1ex}

\paragraph*{\bf Exact verification of binary assignments.}
As these first results show, the values of the binary variables in the \MIP solutions are not supported by the
power vector given by the continuous variables.
In our second experiment, we tried to test whether the binary part of the solutions are correct in the sense that there
exists a power vector $p$ that satisfies these receiver-transmitter assignments.
To this end, we fixed the binary variables to their value in the \MIP solution and solved the remaining \LP, effectively
obtaining a PAP instance as defined in Sec.~\ref{sec:WND}, exactly with \QSOPTEX.  Note that this is a pure feasibility problem.

For the unscaled models, all \LP{}s turned out to be infeasible, as is indicated by the symbol ``$\varnothing$''
in Table~\ref{tab:fixing0}.  On the contrary, the \LP{s} obtained from the scaled models could all be verified as feasible.
Hence the exact \LP solver computed a power vector $p$ to serve all receivers as claimed by the
\MIP~solver.

Additionally, we can see that proving the infeasibility of the unscaled \LP{}s took notably longer than proving the scaled
\LP{}s feasible.  The reason is that in the first case, \QSOPTEX always had to apply increased 128bit arithmetic, while
for the scaled \LP{}s, the basis information after initial double-precision solve turned out to be already exactly feasible.

\vspace*{-1ex}

\paragraph*{\bf Exact \MIP solving.}
We stress that the approach above only yields proven primal bounds on the optimal objective value.  Because \CPLEX uses
floating-point \LP{} bounds, it is unclear whether optimal solutions have been cut off.

In order to further investigate this, we tried to apply the exact extension of the \SCIP solver.  However, for all but
the smallest instances, we could not get any results.  For the instance with 225 transmitters and 20 receivers, the
solving took over 20~hours, 139\,097\,820 branch-and-bound nodes, and more than 7\,GB peak memory usage.
The result was 42 and thus confirmed the optimality of the solution found by \CPLEX.

The slow performance is not really surprising, since the current implementation is a pure branch-and-bound algorithm and
lacks many of the sophisticated features of today's state-of-the-art \MIP solvers.  Hence, this should not be taken as
a proof that exact \MIP solvers are in principle not applicable to this application.

\vspace*{-1ex}

\paragraph*{\bf Accurate computation of the power vector.}
Arguably, computing the power vector exactly is more than necessary for the practical application, and the running times
of \QSOPTEX with almost one minute for the largest \LP may become a bottleneck.  However, in practice it suffices to
compute a power vector that satisfies the original \SIR inequalities~\eqref{eq:firstSIRineq} within a reasonably small tolerance.
In our last experiment, we tested whether the idea of iterative refinement available in the
\SOPLEX solver, can achieve this faster than an exact \LP solver.  We used an (absolute) tolerance of $10^{-25}$, which for
the scaled models suffices to guarantee a tolerance of the same order of magnitude for~\eqref{eq:firstSIRineq}.

Table~\ref{tab:ir} shows the results: the actually reached maximum violation of the \LP rows (as small as
$10^{-40}$), the solving time, and its relative difference to the running times of \QSOPTEX.  For all but the two instances that are solved within one second, \SOPLEX is at least twice as fast as \QSOPTEX.
Note, however, that the implementation of both solvers, in particular the simplex method, differs in many details, and
so we cannot draw a reliable conclusion, let alone on such a limited test set.
However, it suggests that iterative refinement may be more suited to the practical setting of certain applications.

\section{Conclusion}
\label{sec:conclusion}

This paper has tried to highlight a number of numerical issues that must be considered when solving \MIP models
for wireless network design.  We demonstrated that the linearization of the crucial \SIR inequalities in
combination with the definition of feasibility used in floating-point solvers can lead to completely unreliable
results and that an a priori scaling of the constraints can help, but it is not able to make the solutions completely reliable.
We also showed that the current performance of exact \MIP solvers is not sufficient to address the combinatorial difficulty
of these models.
On the positive side, we could show that recent advances in exact and accurate \LP solving are of great help for
computing reliable primal solutions.
So far, we have applied these only as a post processing after the \MIP solution process.  Ideally, however, the
accurate solution of \LP{}s on the continuous variables should be integrated into the branch-and-bound process and used
as a direct verification of the primal bound given by the incumbent solution.
An important next step will be to extend the computational experiments to
a larger set of test instances including other types
of wireless technologies such as DVB-T.


\begin{thebibliography}{10}
\providecommand{\url}[1]{\texttt{#1}}
\providecommand{\urlprefix}{URL }

\bibitem{Achterberg2009}
Achterberg, T.: {SCIP}: {S}olving {C}onstraint {I}nteger {P}rograms. Math.
  Prog. Computation  1(1),  1--41 (2009)

\bibitem{AmEtAl06}
Amaldi, E., Capone, A., Malucelli, F., Mannino, C.: Optimization problems and
  models for planning cellular networks. In: Resende, M., Pardalos, P. (eds.)
  Handbook of Optimization in Telecommunication, pp. 917--939. Springer (2006)

\bibitem{Applegate2007b}
Applegate, D.L., Cook, W., Dash, S., Espinoza, D.G.: {QS}opt\_ex (2007), {\tt
  http://www.dii.uchile.cl/\~{}daespino/ESolver\_doc/}

\bibitem{BaEtAl14}
Bauschert, T., B\"using, C., D'Andreagiovanni, F., Koster, A.M.C.A., Kutschka,
  M., Steglich, U.: Network planning under demand uncertainty with robust
  optimization. IEEE Comm. Mag.  52,  178--185 (2014)

\bibitem{CaEtAl11}
Capone, A., Chen, L., Gualandi, S., Yuan, D.: A new computational approach for
  maximum link activation in wireless networks under the sinr model. IEEE
  Trans. Wireless Comm.  10(5),  1368--1372 (2011)

\bibitem{Cook2011}
Cook, W., Koch, T., Steffy, D.E., Wolter, K.: An exact rational mixed-integer
  programming solver. In: G{\"u}nl{\"u}k, O., Woeginger, G. (eds.) Integer
  Programming and Combinatoral Optimization, LNCS, vol. 6655, pp. 104--116.
  Springer (2011)

\bibitem{CookKochSteffyWolter2013}
Cook, W., Koch, T., Steffy, D.E., Wolter, K.: A hybrid branch-and-bound
  approach for exact rational mixed-integer programming. Math. Prog.
  Computation  (2013)

\bibitem{DA11}
D'Andreagiovanni, F.: On improving the capacity of solving large-scale wireless
  network design problems by genetic algorithms. In: Di~Chio, C., et~al. (eds.)
  EvoApplications 2011. LNCS, vol. 6625, pp. 11--20. Springer, Heidelberg
  (2011)

\bibitem{DA12}
D'Andreagiovanni, F.: Pure 0-1 programming approaches to wireless network
  design. 4OR-Q. J. Oper. Res.  (2012), doi: 10.1007/s10288-011-0162-z

\bibitem{DAMaSa11}
D'Andreagiovanni, F., Mannino, C., Sassano, A.: Negative cycle separation in
  wireless network design. In: Pahl, J., Reiners, T., Voss, T. (eds.) Network
  Optimization, LNCS, vol. 6701, pp. 51--56. Springer, Heidelberg (2011)

\bibitem{DAMaSa13}
D'Andreagiovanni, F., Mannino, C., Sassano, A.: {GUB} covers and power indexed
  formulations for wireless network design. Management Science  59(1),
  142--156 (2013)

\bibitem{Dhiflaoui2003}
Dhiflaoui, M., Funke, S., Kwappik, C., Mehlhorn, K., Seel, M., Sch{\"o}mer, E.,
  Schulte, R., Weber, D.: Certifying and repairing solutions to large {LP}s:
  How good are {LP}-solvers? In: Proc. SODA 2003. pp. 255--256. SIAM (2003)

\bibitem{Espinoza2006}
Espinoza, D.G.: On Linear Programming, Integer Programming and Cutting Planes.
  {Ph.D.} thesis, Georgia Institute of Technology (2006)

\bibitem{GleSteWol2012}
Gleixner, A.M., Steffy, D.E., Wolter, K.: Improving the accuracy of linear
  programming solvers with iterative refinement. In: Proc. ISSAC 2012. Grenoble
  (2012)

\bibitem{GleSteWol2016}
Gleixner, A.M., Steffy, D.E., Wolter, K.: Iterative refinement for linear
  programming. INFORMS J. Computing  (2016), to appear - available as
  ZIB-Report 15-15

\bibitem{gmp}
{GNU} Multiple Precision Arithmetic Library~Version, G.: {\tt
  http://gmplib.org/}

\bibitem{cplex}
{ILOG CPLEX Optimizer}:
  \url{http://www-01.ibm.com/software/integration/optimization/cplex-optimizer/}

\bibitem{KeOlRa10}
Kennington, J., Olinick, E., Rajan, D.: Wireless Network Design: Optimization
  Models and Solution Procedures. Springer, Heidelberg (2010)

\bibitem{Koch2004}
Koch, T.: The final {NETLIB-LP} results. Oper. Res. Lett.  32(2),  138--142
  (2004)

\bibitem{MaRoSm06}
Mannino, C., Rossi, F., Smriglio, S.: The network packing problem in
  terrestrial broadcasting. Oper. Res.  54(6),  611--626 (2006)

\bibitem{MaRoSm07}
Mannino, C., Rossi, F., Smriglio, S.: A unified view in planning broadcasting
  networks. In: DIS Technical Report. vol. 8-2007. Aracne Editrice, Roma (2007)

\bibitem{neumaier2004safe}
Neumaier, A., Shcherbina, O.: Safe bounds in linear and mixed-integer linear
  programming. Mathematical Programming  99(2),  283--296 (2004)

\bibitem{Ra01}
Rappaport, T.S.: Wireless Communications: Principles and Practice. Prentice
  Hall, Upper Saddle River (2001)

\bibitem{scip}
{SCIP}: {Solving Constraint Integer Programs.} \url{http://scip.zib.de}

\bibitem{soplex}
{SoPlex. The Sequential object-oriented simPlex.}: {\tt http://soplex.zib.de/}

\end{thebibliography}

\end{document}